\documentclass[a4paper,12pt]{article}

\usepackage[english]{babel} 
\usepackage[ansinew]{inputenc}
\usepackage[LGR,T1]{fontenc}

\usepackage{verbatim}
\usepackage{color}
\usepackage{calc}
\usepackage{hyperref}
\usepackage[dvips]{graphicx}
\hypersetup{colorlinks=true, linkcolor=blue, filecolor=blue, urlcolor=blue}

\usepackage[mathscr]{eucal}
\usepackage{amscd,amsmath,amsfonts,amssymb}
\everymath{\displaystyle}

\newcommand{\mb}[1]{\mathbb{#1}}
\newcommand{\mc}[1]{\mathcal{#1}}

\newcommand{\MBC}{\mathbb{C}}

\newcommand{\MBQ}{\mathbb{Q}}
\newcommand{\MBR}{\mathbb{R}}

\newcommand{\MBZ}{\mathbb{Z}}

\newcommand{\MCE}{\mathcal{E}}

\newcommand{\MCJ}{\mathcal{J}}

\newcommand{\MCO}{\mathcal{O}}
\newcommand{\MCP}{\mathcal{P}}

\newcommand{\MCR}{\mathcal{R}}
\newcommand{\MCS}{\mathcal{S}}

\newcommand{\MCX}{\mathcal{X}}

\newcommand{\MFa}{\mathfrak{a}}

\newcommand{\MFb}{\mathfrak{b}}

\newcommand{\MFf}{\mathfrak{f}}

\newcommand{\MFm}{\mathfrak{m}}

\newcommand{\MFn}{\mathfrak{n}}

\newcommand{\MFp}{\mathfrak{p}}

\newcommand{\MFq}{\mathfrak{q}}

\newcommand{\MSO}{\mathscr{O}}

\newcommand{\GVe}{\varepsilon}

\newcommand{\GGs}{\sigma}

\newcommand{\GGW}{\Omega}

\newcommand{\vid}{\varnothing}

\newcommand{\lA}{\left\{}

\newcommand{\la}{\left\langle}
\newcommand{\ra}{\right\rangle}

\usepackage{makeidx}
\title{Index-modules and applications.}
\author{St\'ephane VIGUI\'E
\footnote{S.Vigui\'e, Laboratoire de math\'ematiques de Besançon, UMR CNRS 6623, Universit\'e de Franche-Comt\'e, 16 route de Gray, 25030 Besançon cedex, France.
e-mail: \texttt{stephane.viguie@univ-fcomte.fr}}
}
\makeindex

\newtheorem{cro}{cro}[section]
\newtheorem{dfe}{dfe}[section]
\newtheorem{lem}{lem}[section]
\newtheorem{pro}{pro}[section]
\newtheorem{rem}{rem}[section]
\newtheorem{teh}{teh}[section]

\newtheorem{cor}[cro]{Corollary}

\newtheorem{df}[dfe]{Definition}
\newtheorem{lm}[lem]{Lemma}
\newtheorem{pr}[pro]{Proposition}

\newtheorem{rmq}[rem]{Remark}
\newtheorem{theo}[teh]{Theorem}

\numberwithin{equation}{section}

\begin{document}
\maketitle

\noindent{\small\textbf{Abstract.}
Let $K$ be a commutative field, $A\subseteq K$ be a Dedekind ring and $V$ be a $K$-vector space.
For any pair of $A$-lattices $R\neq0$ and $S$ of $V$, we define an $A$-submodule $\left[R:S\right]'_A$ of $K$, their $A$-index-module.
Once the basic properties of these modules are stated, we show that this notion can be used to recover more usual ones: the group-index, the relative invariant, the Fitting ideal of $R/S$ when $S\subseteq R$, and the generalized index of Sinnott.
As an example, we consider the following situation.
Let $F/k$ be a finite abelian extension of global function fields, with Galois group $\mathrm{G}$, and degree $g$.
Let $\infty$ be a place of ${k}$ which splits completely in $F/k$.
Let $\MCO_F$ be the ring of functions of $F$, which are regular outside the places of $F$ sitting over $\infty$.
Then one may use Stark units to define a subgroup $\MCE_F$ of $\MCO_F^\times$, the group of units of $\MCO_F$.
We use the notion of index-module to prove that for every nontrivial irreducible rational character $\psi$ of $\mathrm{G}$, the $\psi$-part of
$\MBZ\left[g^{-1}\right]\otimes_\MBZ\left(\MCO_F^\times/\MCE_F\right)$ and the $\psi$-part of $\MBZ\left[g^{-1}\right]\otimes_\MBZ Cl(\MCO_F)$
have the same order.}
\\


\section{Introduction.}

Let $K$ be a commutative field, $A\subseteq K$ be a Dedekind ring and $V$ be a $K$-vector space.
In this paper, we associate to every pair of $A$-lattices $R\neq0$ and $S$ of $V$ their $A$-index-module
\[[R:S]'_A\,\,:=\,\,\{det(u);u\in End_{K}(V') \quad\text{and}\quad u(R)\subseteq S\},\]
where $V'$ is the $K$-subspace of $V$ generated by $R$ and $S$, $V'=\la R,S\ra_K$.
As we will see below, $[R:S]'_A$ is a finitely generated $A$-submodule of $K$, unless $\la R\ra_K\neq\la S\ra_K$ and $d_R\leq d_S$ (see the notation below).
If $K$ is the fraction field of $A$, and $\la R\ra_K=\la S\ra_K=V$, then $[R:S]'_A$ is equal to the ''relative invariant'' $\chi(S,R)$, defined in \cite[§4, n°6, page 63]{bourbaki.AC7.65}.
Moreover, our $[R:S]'_A$ recover the notion of the ''generalized index'', used for instance in \cite{sinnott80}, \cite{yin97a}, \cite{oukhaba09} and \cite{belliard-nguyen-quang-do05}.
If $S\subseteq R$, then $[R:S]'_A$ is equal to the $A$-Fitting ideal of $R/S$,
\[[R:S]'_A=Fitt_A(R/S).\]
The index-module have nice properties making it very easily be handled.
For example, we consider in section \ref{special units} the following situation.
Let $F/{k}$ be a finite abelian extension of global function fields, with Galois group $\mathrm{G}$, and degree $g$.
Let $\infty$ be a place of ${k}$ which splits completely in $F/{k}$.
Let $\MCO_{k}$ (resp. $\MCO_F$) be the ring of functions of $k$ (resp. $F$), which are regular outside $\infty$ (resp. the places of $F$ sitting over $\infty$).
Then one may use Stark units to define a subgroup $\MCE_F$ of $\MCO_F^\times$, the group of units of $\MCO_F$.
The group $\MCO_F^\times/\MCE_F$ is finite.
Since Stark units are related to $L$-functions, the order of $\MCO_F^\times/\MCE_F$ is related to the order of the ideal class group $Cl(\MCO_F)$ of $\MCO_F$, thanks to the analytic class number formula.
We use the notion of index-module to prove that for every nontrivial irreducible rational character $\psi$ of $\mathrm{G}$, the $\psi$-part of
$\MBZ\left[g^{-1}\right]\otimes_\MBZ\left(\MCO_F^\times/\MCE_F\right)$ and the $\psi$-part of $\MBZ\left[g^{-1}\right]\otimes_\MBZ Cl(\MCO_F)$
have the same order.

\section{Basic properties of the index-modules.}\label{indexmodules}

By an $A$-lattice of $V$, we mean a finitely generated $A$-submodule $R$ of $V$, such that the $K$-vector subspace of $V$ generated by $R$, denoted by $\la R\ra_K$ or $KR$, has dimension equal to the $A$-rank of $R$, $d_R\,\,:=\,\,dim_K\la R\ra_K=rk_A(R)$.
If $R\neq0$, we know that there exists a fractional ideal $\MFm_R$ of $A$, and $B_R=\left(b_{R,0},..., b_{R,d_R-1}\right)$ a $K$-basis of $\la R\ra_K$, such that
$R=\MFm_Rb_{R,0}\oplus\mathop{\oplus}_{i=1}^{d_R-1}Ab_{R,i}$.
Moreover, $\MFm_R$ can be chosen integral.
(See \cite[§4, n°10, Proposition 24, page 79]{bourbaki.AC7.65}.)

In the sequel, $R\neq0$, $S$, and $T$ are $A$-lattices of $V$.
When $A$ is implicit, we will simply note $[R:S]'$.

\begin{rmq}\label{trucidiot}
For any isomorphism $u$ of $V$, we have $\left[u(R):u(S)\right]'_A=[R:S]'_A$.
\end{rmq}

\begin{pr}\label{dfmodinc}
We have the following formulas.
\[[R:S]'_A = \lA\begin{array}{ll}
0 & \text{if} \qquad {d_R}>{d_S} \\
K & \text{if} \qquad {d_R}\leq{d_S} \qquad \text{and} \qquad KS\neq KR \\
{\MFm_S}{\MFm_R}^{-1}det_{B_S,B_R}(Id_{V'}) & \text{if} \qquad V'=KS=KR.
\end{array}\right.\]
\end{pr}

\noindent\textsl{Proof.}
The case $KS\neq KR$ is easy, and left to the reader.
In the sequel, we assume $KS=KR$.
We set $D=det_{B_S,B_R}(Id_{V'})$, and $n={d_R}-1={d_S}-1$ for convenience.

Let $a\in {\MFm_S}{\MFm_R}^{-1}$, and let $u_a$ be the unique automorphism of $V'$ such that $u_a(b_{R,0})=a.b_{S,0}$, and such that for all $i\in\{1,...,n\}$, $u_a({b_{R,i}})={b_{S,i}}$.
Clearly $u_a(R)\subseteq S$, and
\[det_{B_R,B_R}(u_a)=det_{B_R,B_S}(u_a)det_{B_S,B_R}(Id_{V'})=aD.\]
From this we deduce ${\MFm_S}{\MFm_R}^{-1}D\subseteq[R:S]'$.

Conversely, let $u$ be an endomorphism of the vector space $V'$ such that $u(R)\subseteq S$.
Let $M:=mat_{B_R,B_S}(u)$ be the matrix of $u$ with respect to the bases $B_R$ and $B_S$. 
Then, for all $a\in{\MFm_R}$, we have
\[u(ab_{R,0})\in S\quad\Longleftrightarrow\quad\sum_{i=0}^naM_{i,0}{b_{S,i}}\in {\MFm_S} b_{S,0}\oplus\mathop{\oplus}_{i=1}^nA{b_{S,i}}.\]
Hence, we have $M_{0,0}\in {\MFm_S}{\MFm_R}^{-1}$ and $M_{i,0}\in {\MFm_R}^{-1}$ for all $i\in\{1,...,n\}$.
For all $j\in\{1,...,n\}$, we have
\[u(b_{R,j})\in S\quad\Longleftrightarrow\quad\sum_{i=0}^nM_{i,j}{b_{S,i}}\in{\MFm_S} b_{S,0}\oplus\mathop{\oplus}_{i=1}^nA{b_{S,i}}.\]
Hence, we have $M_{0,j}\in {\MFm_S}$ and $M_{i,j}\in A$ for all $i\in\{1,...,n\}$.
Let $\sigma$ be a permutation of $\{0,...,n\}$.
If $\sigma(0)=0$, then for all $i\in\{1,...,n\}$, $\sigma(i)\in\{1,...,n\}$.
In this case, $M_{\sigma(0),0}$ (i.e $M_{0,0}$) belongs to ${\MFm_S}{\MFm_R}^{-1}$, and for all $i\in\{1,...,n\}$, $M_{\sigma(i),i}\in A$.
So $\prod_{i=0}^nM_{\sigma(i),i}\in{\MFm_S}{\MFm_R}^{-1}$.

If $\sigma(0)\in\{1,...,n\}$, then $M_{\sigma(0),0}\in{\MFm_R}^{-1}$.
There is $l\in\{1,...,n\}$ such that $0=\sigma(l)$, and then $M_{\sigma(l),l}$ (i.e $M_{0,l}$) belongs to ${\MFm_S}$.
For all $j\in\{1,...,n\}\setminus\{l\}$, $M_{\sigma(j),j}\in A$.
So $\prod_{i=0}^nM_{\sigma(i),i}\in{\MFm_S}{\MFm_R}^{-1}$.
This proves that $det_{B_R,B_S}(u)\in{\MFm_S}{\MFm_R}^{-1}$.
Since 
\[det(u)=det_{B_R,B_S}(u)det_{B_S,B_R}(Id_{V'})=det_{B_R,B_S}(u)D,\] 
we get $det(u)\in{\MFm_S}{\MFm_R}^{-1}D$.
Thus we have verified that $[R:S]'\subseteq{\MFm_S}{\MFm_R}^{-1}D$, so that $[R:S]'={\MFm_S}{\MFm_R}^{-1}D$.
\hfill $\square$

\begin{rmq}\label{RRRAAegalA}
If $R=S$, then we obtain $[R:R]'_A=A$.
\end{rmq}

\begin{cor}
The $A$-index-module $[R:S]_A'$ is an $A$-module, which is finitely generated unless ${d_R}\leq{d_S}$ and $KS\neq KR$.
\end{cor}

\begin{pr}\label{modincpourreseauxisomorphes}
Assume that $KR=KS$.
Then the following properties are equivalent.

\noindent i) The $A$-modules $R$ and $S$ are isomorphic.

\noindent ii) There is an automorphism $u$ of $V'$ such that $u(R)=S$.

\noindent iii) $[R:S]'$ is a cyclic $A$-module.

\noindent Suppose these properties are satisfied. Then for any automorphism $v$ of $V'$ such that $v(R)=S$, we have $[R:S]'=A.det(v)$.
\end{pr}

\noindent\textsl{Proof.}
$i)\Rightarrow ii)$ is trivial.
If $ii)$ is true, we can choose $(\MFm_R,B_R)$ and $(\MFm_S,B_S)$ such that ${\MFm_R}={\MFm_S}$, and then we deduce $iii)$ from Proposition \ref{dfmodinc}.
Let $x\in K$ be such that $[R:S]'={\MFm_S}{\MFm_R}^{-1}x$.
If $iii)$ is true, then ${\MFm_S}{\MFm_R}^{-1}$ is principal, so ${\MFm_R}$ and ${\MFm_S}$ are isomorphic as $A$-modules, and $i)$ follows.

Finally, assume $i)$, $ii)$, and $iii)$ are true.
Let $v$ be an automorphism of $V'$ such that $v(R)=S$.
We choose $(\MFm_R,B_R)$ and $(\MFm_S,B_S)$ such that ${\MFm_S}={\MFm_R}$, and ${b_{S,i}}=v({b_{R,i}})$ for all $i\in\{0,...,{d_R}-1\}$.
Then Proposition \ref{dfmodinc} gives 
\[[R:S]'=A.det_{B_S,B_R}(Id_{V'})=A.det_{B_R,B_S}(v)det_{B_S,B_R}(Id_{V'})=A.det(v).\]
\hfill $\square$

\begin{rmq}
In \cite{sinnott80}, W.~Sinnott defined a generalized index $(R:S)$, for the case $d_R=d_S=dim_K(V)<\infty$, when $A=\MBZ$ and $K=\MBQ$ or $K=\MBR$, or when $A=\MBZ_p$ and $K=\MBQ_p$, with $p$ a prime number.
By definition, if $u$ is an endomorphism of $V$ such that $u(R)=S$, then:

- $(R:S)=|det(u)|$ in the case $A=\MBZ$ and $K=\MBQ$ or $K=\MBR$.

- $(R:S)=p^{v\left(det(u)\right)}$ in the case $A=\MBZ_p$ and $K=\MBQ_p$, with $v$ the normalized valuation on $\MBQ_p$.

From Proposition \ref{modincpourreseauxisomorphes}, we have $[R:S]'_A=A(R:S)$ in both cases.
\end{rmq}

\begin{pr}\label{prop calcul modinc}
Let $n=dim_K(V)$, $(\MFm_i)_{i=1}^n$ and $(\MFn_i)_{i=1}^n$ two families of fractional ideals of $A$, and $B=(b_i)_{i=1}^n$ a $K$-basis of $V$.
We set $M=\mathop{\oplus}_{i=1}^n\MFm_ib_i$ and $N=\mathop{\oplus}_{i=1}^n\MFn_ib_i$.
Then 
\[[M:N]'=\prod_{i=1}^n\MFn_i\MFm_i^{-1}.\]
\end{pr}

\noindent\textsl{Proof.}
For all $i\in\{1,...,n\}$, let $x_i\in \MFn_i\MFm_i^{-1}$.
Let $u$ be the endomorphism of $V$ such that for all $i\in\{1,...,n\}$, $u(b_i)=x_ib_i$.
Then $u(M)\subseteq N$, and $det(u)=\prod_{i=1}^nx_i$.
This shows that $\prod_{i=1}^n\MFn_i\MFm_i^{-1}\subseteq[M:N]'$.

Let $v$ be an endomorphism of $V$ such that $v(M)\subseteq N$.
Let $P=mat_{B,B}(v)$.
For all $j\in\{1,...,n\}$, and all $a\in\MFm_j$, we have
\[v(ab_j)\in S\quad\Longleftrightarrow\quad\sum_{i=1}^nP_{i,j}ab_i\in\mathop{\oplus}_{i=1}^n\MFn_ib_i.\]
So $P_{i,j}\in\MFn_i\MFm_j^{-1}$, for all $i\in\{1,...,n\}$.
For all permutation $\sigma$ of $\{1,...,n\}$, we have $\prod_{i=1}^nP_{\sigma(i),i}\in\prod_{i=1}^n\MFn_{\sigma(i)}\MFm_i^{-1}$, and so $\prod_{i=1}^nP_{\sigma(i),i}\in\prod_{i=1}^n\MFn_{i}\MFm_i^{-1}$.
Then $det(v)\in\prod_{i=1}^n\MFn_{i}\MFm_i^{-1}$.
Thus, we have verified $[R:S]'\subseteq\prod_{i=1}^n\MFn_{i}\MFm_i^{-1}$.
\hfill $\square$

\begin{pr}\label{propprodmodinc}
Assume $KR=KS$ or $KS=KT\neq0$.
Then
\[[R:T]'=[R:S]'[S:T]'.\] 
(This product has to be understood in the following way:
$[R:T]'$ is the sub-$A$-module of $K$ generated by all the elements $xy$, with $(x,y)\in[R:S]'\times[S:T]'$.)
\end{pr}

\noindent\textsl{Proof.}
The inclusion $[R:S]'[S:T]'\subseteq[R:T]'$ is trivial.
If $[R:S]'$ is zero, then ${d_S}<{d_R}$, and $KS=KT$.
So ${d_T}<{d_R}$, and $[R:T]'$ is zero.
A similar argument shows that if $[S:T]'$ is zero, then $[R:T]'$ is zero.

If $[R:S]'\neq0$ with $KR\neq KS$, then $KS=KT$, $[R:S]'=K$ and $[S:T]'\neq0$.
This implies $K=[R:S]'[S:T]'\subseteq[R:T]'$, and thus $[R:T]'=K$.
In the same way, if $[S:T]'\neq0$ with $KS\neq KT$, then $[R:T]'=K$.

To conclude, it suffices to study the case $KR=KS=KT$.
Then, by Proposition \ref{dfmodinc}, we have
\[\begin{array}{lllll}
[R:T]' & = & {\MFm_T}{\MFm_R}^{-1}det_{B_T,B_R}(Id_{V'}) & & \\
 & = & {\MFm_S}{\MFm_R}^{-1}det_{B_S,B_R}(Id_{V'}){\MFm_T}{\MFm_S}^{-1}det_{B_T,B_S}(Id_{V'}) & = & [R:S]'[S:T]'.
\end{array}\]
\hfill $\square$

\begin{cor}
Assume that $KR=KS$.
Then $[R:S]'=\left([S:R]'\right)^{-1}$, where $\left([S:R]'\right)^{-1}$ is the sub-$A$-module $\left\{x\in K;\,[R:S]'x\subseteq A\right\}$ of $K$.
\end{cor}

\begin{lm}\label{lemmadeuxpointdouze}
Let $V''$ be a nonzero sub-$K$-vector space of $V$.
Let $n=dim_K(V'')$, $(\MFm_i)_{i=1}^n$ and $(\MFn_i)_{i=1}^n$ two families of fractional ideals of $A$, $B=(b_i)_{i=1}^n$ and $C=(c_i)_{i=1}^n$ two $K$-basis of $V''$.
We set $M=\mathop{\oplus}_{i=1}^n\MFm_ib_i$ and $N=\mathop{\oplus}_{i=1}^n\MFn_ic_i$.
Then 
\[[M:N]'=\left(\prod_{i=1}^n\MFn_i\MFm_i^{-1}\right)det_{C,B}(Id_{V''}).\]
\end{lm}

\noindent\textsl{Proof.}
By Proposition \ref{propprodmodinc}, we have
\[[M:N]'=[M:\mathop{\oplus}_{i=1}^nAb_i]'[\mathop{\oplus}_{i=1}^nAb_i:\mathop{\oplus}_{i=1}^nAc_i]'[\mathop{\oplus}_{i=1}^nAc_i:N]'.\]
From Proposition \ref{prop calcul modinc}, we obtain $[M:\mathop{\oplus}_{i=1}^nAb_i]'=\prod_{i=1}^n\MFm_i^{-1}$, and $[\mathop{\oplus}_{i=1}^nAc_i:N]'=\prod_{i=1}^n\MFn_i$.
By Proposition \ref{modincpourreseauxisomorphes}, $[\mathop{\oplus}_{i=1}^nAb_i:\mathop{\oplus}_{i=1}^nAc_i]'=A.det(u)$, where $u$ is the unique endomorphism of $V''$ such that for all $i\in\{1,...,n\}$, $u(b_i)=c_i$.
But $det(u)=det_{C,B}(Id_{V''})$, and the lemma follows.
\hfill $\square$

\begin{pr}\label{prdirectsums}
(Direct sums)
Let $n\in\mb{N}^\ast$, $(R_i)_{i=1}^n$ a family of sub-$A$-modules of $R$, and $(S_i)_{i=1}^n$ a family of sub-$A$-modules of $S$, such that
\[R=\mathop{\oplus}_{i=1}^nR_i\quad\text{and}\quad S=\mathop{\oplus}_{i=1}^nS_i.\]
We assume that for all $i\in\{1,...,n\}$, $KS_i\subseteq KR_i\neq0$.
Then
\[[R:S]'=\prod_{i=1}^n[R_i:S_i]'.\]
\end{pr}

\noindent\textsl{Proof.}
We have $KR=\mathop{\oplus}_{i=1}^nKR_i$ and $KS=\mathop{\oplus}_{i=1}^nKS_i$.
If there is $i\in\{1,...,n\}$ such that $KS_i\neq KR_i$, then $KR\neq KS$, $[R:S]'=0$ and $[R_i:S_i]'=0$, so the result is trivial.
Suppose now that for all $i\in\{1,...,n\}$, $KS_i=KR_i$.
For all $i\in\{1,...,n\}$, let $q_i=dim_K(KR_i)$, $B_i=(b_{i,j})_{j=1}^{q_i}$ and $C_i=(c_{i,j})_{j=1}^{q_i}$ two $K$-basis of $KR_i$, $(\MFm_{i,j})_{j=1}^{q_i}$ and $(\MFn_{i,j})_{j=1}^{q_i}$ two families of fractional ideals of $A$ such that
\[R_i=\mathop{\oplus}_{j=1}^{q_i}\MFm_{i,j}b_{i,j}\quad\text{and}\quad S_i=\mathop{\oplus}_{j=1}^{q_i}\MFn_{i,j}c_{i,j}.\]
We also set $B=(b_{1,1},...,b_{1,q_1},...,b_{n,1},...,b_{n,q_n})$ and $C=(c_{1,1},...,c_{1,q_1},...,c_{n,1},...,c_{n,q_n})$.
These are two $K$-basis of $KR$.
Then we have $[R_i:S_i]'=\left(\prod_{j=1}^{q_i}\MFn_{i,j}\MFm_{i,j}^{-1}\right)det_{C_i,B_i}(Id_{V_i})$, thanks to Lemma \ref{lemmadeuxpointdouze}.
But we have the decompositions
\[R=\mathop{\oplus}_{i=1}^n\mathop{\oplus}_{j=1}^{q_i}\MFm_{i,j}b_{i,j}\quad\text{and}\quad S=\mathop{\oplus}_{i=1}^n\mathop{\oplus}_{j=1}^{q_i}\MFn_{i,j}c_{i,j}.\]
Hence, from Lemma \ref{lemmadeuxpointdouze} we deduce
\[\begin{array}{lllll}
[R:S]' & = & det_{C,B}(Id_{V'})\prod_{i=1}^n\prod_{j=1}^{q_i}\MFn_{i,j}\MFm_{i,j}^{-1} & & \\
 & = & \prod_{i=1}^n\left(\prod_{j=1}^{q_i}\MFn_{i,j}\MFm_{i,j}^{-1}\right)det_{C_i,B_i}(Id_{V_i}) & = & \prod_{i=1}^n[R_i:S_i]'.
\end{array}\]
\hfill $\square$

\begin{pr}\label{propextscal pour modinc}
(Scalar expansion)
Let $B$ be a Dedekind ring, embedded in $K$, such that $A$ is a sub-ring of $B$.
Then $[BR:BS]'_B=B[R:S]'_A$.
(Where $B[R:S]'_A$ is the $B$-submodule of $K$ generated by $[R:S]'_A$.)
\end{pr}

\noindent\textsl{Proof.}
First, notice that $BR$ is a $B$-lattice of $V$, thanks to the inequalities
\[rk_A(R) = dim_K(KR) \leq rk_B(BR) \leq rk_A(R).\]
If ${d_S}<{d_R}$, then $[BR:BS]'_B=0=[R:S]'_A$.
If ${d_R}\leq {d_S}$, with $KS\neq KR$, then $[BR:BS]'_B=K=[R:S]'_A$.
Suppose now $KR=KS$.
By Proposition \ref{dfmodinc}, we have
\[[BR:BS]'_B=(B{\MFm_S})(B{\MFm_R})^{-1}det_{B_S,B_R}(Id_{V'})=B.{\MFm_S}{\MFm_R}^{-1}det_{B_S,B_R}(Id_{V'})=B[R:S]'_A.\]
\hfill $\square$

\begin{pr}\label{linkindmodinc}
Assume $S\subseteq R$ and $d_R=d_S$.
If $A/\MFa$ is finite for any nonzero ideal $\MFa$ of $A$, then
\[\left[A:[R:S]'_A\right]=[R:S],\]
where $\left[A:[R:S]'_A\right]$ (resp. $[R:S]$) is the group-index of the ideal $[R:S]'_A$ in $A$ (resp. of $S$ in $R$).
\end{pr}

\noindent\textsl{Proof.}
If $A$ is a field then $R=S$ and $[R:S]'_A=A$.
In the sequel, we assume $A$ is not a field, and we consider two cases.

First, suppose $A$ is a discrete valuation ring.
Let $\pi$ be a uniformizer of $A$.
We can choose $\MFm_R=A$, and $B_R$ such that there is $(m_i)_{i=1}^{d_R}\in\mb{N}^{d_R}$, with $S=\mathop{\oplus}_{i=1}^{d_R}\pi^{m_i}b_{R,i-1}$. Then $R/S\simeq\mathop{\oplus}_{i=1}^nA/\pi^{m_i}A$, and $[R:S]'_A=\pi^{\sum_{i=1}^nm_i}A$.
This implies the proposition if $A$ is a discrete valuation ring.

Suppose now $A$ is any Dedekind ring.
We have $R/S\simeq\mathop{\oplus}_{\MFp\in Spec^\ast(A)}A_\MFp\otimes_A(R/S)$ (where $Spec^\ast(A)$ is the set of nonzero prime ideals of $A$).
Since $A_\MFp$ is a flat $A$-module and a discrete valuation ring, for all $\MFp\in Spec^\ast(A)$, we deduce
\[\begin{array}{lllll}
[R:S] & = & \prod_{\MFp\in Spec^\ast(A)}[A_\MFp R:A_\MFp S] & = & \prod_{\MFp\in Spec^\ast(A)}\left[A_\MFp:[A_\MFp R:A_\MFp S]'_{A_\MFp}\right].
\end{array}\]
Taking into account the proposition \ref{propextscal pour modinc}, we obtain
\[\begin{array}{lllll}
\prod_{\MFp\in Spec^\ast(A)}\left[A_\MFp:[A_\MFp R:A_\MFp S]'_{A_\MFp}\right] & = & \prod_{\MFp\in Spec^\ast(A)}\left[A_\MFp:A_\MFp[R:S]'_A\right] & = & \left[A/[R:S]'_A\right].
\end{array}\]
\hfill $\square$

\begin{cor}\label{indmodinc}
Suppose $A=\MCS^{-1}\MBZ$, where $\MCS$ is a multiplicative part of the ring $\MBZ$, which does not contain $0$.
If $S\subseteq R$ and $d_R=d_S$, we have $[R:S]'_A=[R:S]A$.
\end{cor}

\noindent\textsl{Proof.}
For any nonzero ideal $\MFa$ of $A$, we have $\MFa=[A:\MFa]A$ and $A/\MFa$ is finite.
From Proposition \ref{linkindmodinc}, we deduce
\[[R:S]'_A = \left[A:[R:S]'_A\right]A = [R:S]A.\]
\hfill $\square$

\bigskip

For the rest of this section, assume $K$ is the fraction-field of $A$, and $KR=KS=V$ (in this case, $[R:S]'$ is a fractional ideal of $A$).
Set $n=dim_K(V)-1$.
For any choice of a $K$-basis of $V$, we have a canonical isomorphism $\bigwedge^{n+1}V\simeq K$.
Through this isomorphism, the canonical image of $R^{n+1}$ in $\bigwedge^{n+1}V$ is identified to a fractional ideal $\MFa$ of $A$.
In the same way, $S^{n+1}$ defines a fractional ideal $\MFb$ of $A$.
Then we define the relative invariant $\chi(S,R)$ of $S$ and $R$, $\chi(S,R) := \MFb\MFa^{-1}$.
It does not depend on the choice of the $K$-basis of $V$ (see \cite[§4, n°6]{bourbaki.AC7.65}, for more details and basic properties about the relative invariant).
We use a multiplicative notation for the relative invariant, instead of the usual additive one, because it is more adapted to our situation.

\begin{pr}
$\chi(S,R)=[R:S]'$. 
\end{pr}

\noindent\textsl{Proof.}
Set $M=\mathop{\oplus}_{i=0}^nA{b_{R,i}}$, and $N=\mathop{\oplus}_{i=0}^nA{b_{S,i}}$.
$b_{R,0}\wedge\dots\wedge b_{R,n}$ is a basis of $\bigwedge^{n+1}V$.
The $A$-lattice of $\bigwedge^{n+1}V$, generated by the canonical image of $M^{n+1}$, is $Ab_{R,0}\wedge\dots\wedge Ab_{R,n}$, and the $A$-lattice of $\bigwedge^{n+1}V$, generated by the canonical image of $R^{n+1}$ is ${\MFm_R} b_{R,0}\wedge\dots\wedge b_{R,n}$.
By definition, $\chi(M,R)={\MFm_R}^{-1}$.
In the same way, $\chi(S,N)={\MFm_S}$.

Let $u$ be the unique endomorphism of $V$ such that $u({b_{R,i}})={b_{S,i}}$ for all $i\in\{0,\dots,n\}$.
Then, applying \cite[§4, n°6, Proposition 13]{bourbaki.AC7.65}: 
\[\chi(N,M)=\left(det(u)\right)=\left(det_{B_S,B_R}(Id_V)det_{B_R,B_S}(u)\right)=\left(det_{B_S,B_R}(Id_V)\right)\]
Finally
\[\chi(S,R)=\chi(S,N)\chi(N,M)\chi(M,R)={\MFm_S}{\MFm_R}^{-1}det_{B_S,B_R}(Id_V)=[R:S]'.\]
\hfill $\square$

\begin{lm}\label{proplienmodincfitt}
Let $0\rightarrow N \rightarrow L \rightarrow M \rightarrow 0$ be an exact sequence of finitely generated $A$-modules, and assume $L$ is free.
$L$ and $N$ are viewed as $A$-lattices of $K\otimes_AL$.
Then
\[Fitt_A(M)=[L:N]'_A.\]
\end{lm}

\noindent\textsl{Proof.}
Let $C=\{c_1,\dots, c_n\}$ be an $A$-basis of $L$.
Let $\mc{M}$ be the set of square $n\times n$-matrices $T:=(t_{i,j})$, with coefficients in $A$, such that for all $j\in\{1,...,n\}$, $\left(\sum_{i=1}^nt_{i,j}c_i\right)\in N$.
Then $Fitt_A(M)$ is the ideal of $A$ generated by the determinants of the matrices $T\in\mc{M}$.
Let $T\in\mc{M}$.
There is a unique endomorphism $f$ of $K\otimes_AL$, such that for all $j\in\{1,...,n\}$, $f(c_j)=\sum_{i=1}^nt_{i,j}c_i$.
Then $f(L)\subseteq N$, and $det(f)=det(T)$.
We deduce the inclusion $Fitt_A(M)\subseteq[L:N]'_A$.

Let $f\in End_{K}(K\otimes_AL)$, such that $f(L)\subseteq N$.
Let $T$ be the matrix of $f$ in the $K$-basis $C$ of $K\otimes_AL$.
Obviously, $T\in\mc{M}$, and so $det(f)\in Fitt_A(M)$.
Thus we have proved $[L:N]'_A\subseteq Fitt_A(M)$.
\hfill $\square$

\begin{theo}\label{theolienmodincfitt}
Let $M$ be a nonzero finitely generated $A$-module, torsion-free over $A$, and $N$ an $A$-submodule of $M$.
$M$ and $N$ are viewed as $A$-lattices of $K\otimes_AM$.
Then $Fitt_A(M/N)=[M:N]'_A$.
\end{theo}

\noindent\textsl{Proof.}
Let us write $M=\MFm b_0\oplus\mathop{\oplus}_{i=1}^nAb_i$, where $\MFm$ is a nonzero integral ideal of $A$ and $(b_0,...,b_n)$ is a $K$-basis of $K\otimes_AM$.
We set $L=\mathop{\oplus}_{i=0}^nAb_i$.
Since $A$ is a Dedekind ring, the exact sequence 
\[0 \rightarrow M/N \rightarrow L/N \rightarrow L/M \rightarrow 0\] 
gives $Fitt_A(L/M)Fitt_A(M/N)= Fitt_A(L/N)$.
By Lemma \ref{proplienmodincfitt} we have 
\[[L:M]'_AFitt_A(M/N) = [L:N]'_A.\]
But $KM=KL$.
Thus, multiplying by $[M:L]'_A$ and applying Proposition \ref{propprodmodinc} and Remark \ref{RRRAAegalA}, we obtain the desired formula.
\hfill $\square$

\section{The index of Stark units in function fields.}\label{special units}

In this section, we apply the notion of index-module to prove Theorem \ref{theocalcind} below, which gives a weak form of the Gras conjecture in positive characteristic.
For this we shall use the following notation.

Let ${k}$ be a global function field, with field of constants $\mb{F}_q$.
Let $\infty$ be a place of ${k}$, of degree $d$ over $\mb{F}_q$.
Then, we denote by $k_\infty$ the completion of $k$ at $\infty$.
Let $\mc{O}_{k}$ be the Dedekind ring of functions $f\in k$ regular outside $\infty$.
Let us also fix $F\subseteq k_\infty$, a finite abelian extension of $k$ in which $\infty$ splits completely, with Galois group $\mathrm{G}$ and degree $g$.

For any finite abelian extension $K$ of ${k}$, we denote by $\mc{O}_K$ the integral closure of $\mc{O}_{k}$ in $K$, and by $\MCO_K^\times$ the group of units of $\mc{O}_K$.
We denote by ${\mu}(K)$ the group of roots of unity in $K$, and by $Cl\left(\MCO_K\right)$ the ideal class group of $\MCO_K$.

\subsection{Stark units in function fields.}\label{elliptic units}

If $\MFm$ is a nonzero ideal of $\MCO_k$ then we denote by $H_\MFm\subseteq k_\infty$ the maximal abelian extension of $k$ contained in $k_\infty$, such that the conductor of $H_\MFm/k$ divides $\MFm$.
In particular, $\infty$ splits completely in $H_\MFm/k$.
The function field version of the abelian conjectures of Stark, proved by P.~Deligne in \cite{tate84} by using \'etale cohomology or by D.~Hayes in \cite{hayes85} by using Drinfel'd modules, claims that, for any proper nonzero ideal $\MFm$ of $\MCO_k$, there exists an element $\GVe_\MFm\in H_\MFm$, unique up to roots of unity, such that

\noindent (i) The extension $H_\MFm\left(\GVe_\MFm^{1/w_\infty}\right)/k$ is abelian, where $w_\infty:=q^d-1$.

\noindent (ii) If $\MFm$ is divisible by two prime ideals then $\GVe_\MFm$ is a unit of $\MCO_{H_\MFm}$.
If $\MFm=\MFq^e$, where $\MFq$ is a prime ideal of $\MCO_k$ and $e$ is a positive integer, then
\[\GVe_\MFm\MCO_{H_\MFm}=(\MFq)_\MFm^{\frac{w_\infty}{w_k}},\]
where $w_k:=q-1$ and $(\MFq)_\MFm$ is the product of the prime ideals of $\MCO_{H_\MFm}$ which divide $\MFq$.

\noindent (iii) We have
\begin{equation} 
L_\MFm(0,\chi)=\frac{1}{w_\infty}\sum_{\GGs\in \mathrm{Gal}(H_\MFm/k)}\chi(\GGs)v_\infty\left(\GVe_\MFm^\GGs\right) 
\label{equaLMFmGVEMFm} 
\end{equation}
for all complex irreducible characters of $\mathrm{Gal}\left(H_\MFm/k\right)$, where $v_\infty$ is the normalized valuation of $k_\infty$.

Let us recall that $s\mapsto L_\MFm(s,\chi)$ is the $L$-function associated to $\chi$, defined for the complex numbers $s$ such that $Re(s)>1$ by the Euler product
\[L_\MFm(s,\chi)=\prod_{v\nmid\MFm}\left(1-\chi(\GGs_v)N(v)^{-s}\right)^{-1},\]
where $v$ describes the set of places of $k$ not dividing $\MFm$.
For such a place, $\GGs_v$ and $N(v)$ are the Frobenius automorphism of $H_\MFm/k$ and the order of the residue field at $v$ respectively.
Let us remark that $\GGs_\infty=1$ and $N(\infty)=q^d$.

For any finite abelian extension $L$ of $k$ we denote by $\MCJ_L\subseteq\mb{Z}\left[\mathrm{Gal}(L/k)\right]$ the annihilator of ${\mu}(L)$.
The description of $\MCJ_L$ given in \cite[Lemma 2.5]{hayes85} and the property (i) of $\GVe_\MFm$ implies that for any $\eta\in\MCJ_{H_\MFm}$ there exists $\GVe_\MFm(\eta)\in H_\MFm$ such that
\[\GVe_\MFm(\eta)^{w_\infty}=\GVe_\MFm^\eta.\]

\begin{df}
Let ${\MCP_F}$ be the subgroup of $F^\ast$ generated by ${\mu}(F)$ and by all the norms 
\[N_{H_\MFm/H_\MFm\cap F}\left(\GVe_\MFm(\eta)\right),\] 
where $\MFm$ is any nonzero proper ideal of $\MCO_{k}$, and $\eta\in\MCJ_{H_\MFm}$.
By definition, the group of Stark units is
\[{\MCE_F}={\MCP_F}\cap \MCO_F^\times.\]
\end{df}

\begin{rmq}
The index $\left[\MCO_F^\times:\MCE_F\right]$ is finite.
This will be proved in the next subsection.
In \cite{oukhaba92}, H.~Oukhaba succeeded in computing this index in case $F\subseteq H_{(1)}$.
He obtained the following formula
\[\left[\MCO_F^\times:\MCE_F\right]=\frac{h\left(\MCO_F\right)}{\left[H_{(1)}:F\right]},\]
where $h\left(\MCO_F\right)$ is the ideal class number of $\MCO_F$.

Let $S$ be a set of places of $k$, which contains $\infty$.
In \cite{popescu99b}, C.~Popescu defined a group $\MCE_S$ of $S$-units of $F$ by using Rubin-Stark units.
If $S=\{\infty\}$, he proved that for any prime number $p$, $p\nmid g$, and every nontrivial irreducible $p$-adic character $\psi$ of $\mathrm{G}$, the Gras conjecture is verified, 
\[\#\left(\MBZ_p\otimes_\MBZ\left(\MCO_F^\times/\MCE_S\right)\right)_\psi = \#\left(\MBZ_p\otimes_\MBZ Cl\left(\MCO_F\right)\right)_\psi,\]
where the subscript $\psi$ means we take the $\psi$-parts.
See \cite[Theorem 3.10]{popescu99b}.

In the sequel, we use index-modules to prove a weak form of the analoguous statement for the group $\MCE_F$ (see Theorem \ref{theocalcind}), i.e for rational characters.
It can be shown that $\MBZ\left[g^{-1}\right]\otimes_\MBZ\MCE_F$ is included in $\MBZ\left[g^{-1}\right]\otimes_\MBZ\MCE_S$.
From Theorem \ref{theocalcind}, it follows that this inclusion is an equality, so that the full Gras conjecture is also true for $\MCE_F$.
\end{rmq}

\subsection{An index formula for Stark units.}\label{indexformulaelliptic}

Let $\ell_F:F^\times\rightarrow\MBZ[\mathrm{G}]$ be the $\mathrm{G}$-equivariant map defined by 
\[\ell_F(x)=\sum_{\GGs\in G}v_\infty\left(x^\GGs\right)\GGs^{-1}.\]
Let $\mu_g$ be the group of $g$-th roots of unity in the field of complex numbers.
Let $\MSO$ be the integral closure of the principal ring $\MBZ_{\la g\ra}:=\MBZ[g^{-1}]$ in $\mb{Q}(\mu_g)$.
Let us denote by $\widehat{\mathrm{G}}$ the group of complex irreducible characters of $\mathrm{G}$.
Then, for every $\chi\in\widehat{\mathrm{G}}$ the idempotent $e_\chi:=\frac{1}{g}\sum_{\GGs\in G}\chi(\GGs)\GGs^{-1}$ belongs to $\MSO[\mathrm{G}]$.
Moreover, if $\zeta\in\mu_g$ is such that $\zeta\neq1$, then $1-\zeta\in\MSO^\times$, thanks to the formula $g=\prod_{\substack{\zeta\in\mu_g\\\zeta\neq1}}(1-\zeta)$.

\begin{df}
Let $\GGW$ be the $\MBZ[\mathrm{G}]$-submodule of $F^\ast$ generated by $\mu(F)$ and the elements of the form $N_{H_\MFm/F\cap H_\MFm}(\GVe_\MFm)$, where $\MFm$ is any nonzero proper ideal of $\MCO_k$.
\end{df}

\begin{pr}\label{chicpsteomega}
Let $\chi\in\widehat{\mathrm{G}}$ be such that $\chi\neq1$.
Let $\chi_{pr}$ be the character of $\mathrm{Gal}\left(H_{\MFf_\chi}/k\right)$ deduced from $\chi$, where $\MFf_\chi$ is the conductor of the fixed field $F_\chi$ of $\mathrm{Ker}(\chi)$.
Then
\begin{equation}
\MSO\ell_F(\GGW)e_\chi=\MSO w_\infty L_{\MFf_\chi}(0,\bar{\chi}_{pr})e_\chi,
\label{equaMSOellGGWechi}
\end{equation}
where $\MSO\ell_F(\GGW)\subseteq\MSO[\mathrm{G}]$ is the $\MSO$-module generated by $\ell_F(\GGW)$.
\end{pr}

\noindent\textsl{Proof.}
The equality (\ref{equaMSOellGGWechi}) is a direct consequence of the property (iii) of Stark units stated above.
We take our inspiration from the computation made in \cite{oukhaba09}.
Let $\MFm$ be a nonzero proper ideal of $\MCO_{k}$ and let $\GVe_{F,\MFm}:=N_{H_\MFm/F\cap H_\MFm}(\GVe_\MFm)$.
If $\chi$ is not trivial on $\mathrm{Gal}\left(F/F\cap H_\MFm\right)$, then
\[\ell_F(\GVe_{F,\MFm})e_\chi=0.\]
But, if $\chi$ is trivial on $\mathrm{Gal}\left(F/F\cap H_\MFm\right)$ then $F_\chi\subseteq H_\MFm$, $\MFf_\chi|\MFm$ and 
\[\ell_F\left(\GVe_{F,\MFm}\right)e_\chi = \left[F:F\cap H_\MFm\right] \left(\sum_{\GGs\in \mathrm{Gal}\left(H_\MFm/k\right)}\overline{\chi'(\GGs)}v_\infty\left(\GVe_{F,\MFm}^\GGs\right)\right)e_\chi,\]
where $\chi'$ is the complex character of $\mathrm{Gal}(H_\MFm/k)$ deduced from $\chi$.
Therefore the equality (\ref{equaLMFmGVEMFm}) gives
\[\ell_F\left(\GVe_{F,\MFm}\right)e_\chi = w_\infty \left[F:F\cap H_\MFm\right] L_\MFm(0,\overline{\chi'})e_\chi.\]
Now, the relation
\[L_\MFm(0,\overline{\chi}) = \prod_{\substack{v|\MFm\\v\nmid\MFf_\chi}}\left(1-\overline{\chi}_{pr}(\GGs_v)\right) L_{\MFf_\chi}(0,\overline{\chi}_{pr})\]
clearly shows that $\ell_F(\GVe_{F,\MFm})e_\chi\in \MSO w_\infty L_{\MFf_\chi}(0,\overline{\chi}_{pr})e_\chi$.
Conversely, the hypothesis that $\chi$ is not trivial implies that there exists some prime ideal $\MFp$ of $\MCO_k$ such that $\chi_{pr}(\GGs_\MFp)\neq1$.
Thus, if we put $\MFm:=\MFp\MFf_\chi$ then
\[\ell_F\left(\GVe_{F,\MFm}\right)e_\chi = w_\infty \left[F:F\cap H_\MFm\right] \left(1-\overline{\chi}_{pr}(\GGs_\MFp)\right) L_{\MFf_\chi}(0,\overline{\chi}_{pr})e_\chi.\]
But, since $[F:F\cap H_\MFm]$ and $\left(1-\overline{\chi}_{pr}(\GGs_\MFp)\right)$ are in $\MSO^\times$ we obtain $w_\infty L_{\MFf_\chi}(0,\overline{\chi}_{pr})e_\chi\in\MSO\ell_F\left(\GGW\right)e_\chi$.
The proposition is now proved.
\hfill $\square$

\bigskip

If $\psi$ is an irreducible character of $\mathrm{G}$, then we denote by $\MCX_\psi$ the set of $\chi\in\widehat{\mathrm{G}}$ such that $\chi|\psi$.
Let $M$ be a $\MBZ_{\la g\ra}[\mathrm{G}]$-module.
Thus we put $M_\psi=e_\psi M$, where $e_\psi=\sum_{\chi\in\MCX_\psi}e_\chi$.
If $M$ is an $\MSO[\mathrm{G}]$-module and $\chi\in\widehat{\mathrm{G}}$ then we put $M_\chi:=e_\chi M$.

\begin{cor}\label{corOMegaJFOmega}
Let $\psi\neq1$ be an irreducible rational character of $\mathrm{G}$.
Then $\left(\MSO\ell_F(\GGW)\right)_\psi$ and $\left(\MSO\MCJ_F\ell_F(\GGW)\right)_\psi$ are $\MSO$-lattices of the $\MBC$-vector space $\MBC[\mathrm{G}]$.
Moreover,
\[\left[\left(\MSO\ell_F(\GGW)\right)_\psi : \left(\MSO\MCJ_F\ell_F(\GGW)\right)_\psi\right]'_\MSO = \MSO\#\left(\MBZ_{\la g\ra} \otimes_\MBZ \mu(F)\right)_\psi.\]
\end{cor}

\noindent\textsl{Proof.}
On one hand, we have the decomposition
\begin{equation}
\left(\MSO\ell_F(\GGW)\right)_\psi = \mathop{\oplus}_{\chi\in\MCX_\psi} \left(\MSO\ell_F(\GGW)\right)_\chi = \mathop{\oplus}_{\chi\in\MCX_\psi} \left(\MSO w_\infty L_{\MFf_\chi}\left(0,\bar{\chi}_{pr}\right)\right)_\chi,
\label{equadeclatttiecs}
\end{equation}
the last equality being an application of Proposition \ref{chicpsteomega}.
On the other hand, since $\psi\neq1$ the primitive character $\chi_{pr}$ is nontrivial for all $\chi\in\MCX_\psi$.
As a consequence, $L_{\MFf_\chi}(0,\bar{\chi}_{pr})\neq0$.
This implies that $\left(\MSO\ell_F(\GGW)\right)_\psi$ is a free $\MSO$-module of rank $\#(\MCX_\psi)=dim(\psi)$, and hence, it is an $\MSO$-lattice of $\MBC[\mathrm{G}]$.
Similar arguments may be used to verify that $\left(\MSO\MCJ_F\ell_F\left(\GGW\right)\right)_\psi$ is an $\MSO$-lattice of $\MBC[\mathrm{G}]$.
Furthermore, using Proposition \ref{prdirectsums}, Proposition \ref{chicpsteomega} and Remark \ref{trucidiot}, we obtain
\begin{eqnarray}
\left[ \left(\MSO \ell_F\left(\GGW\right)\right)_\psi : \left(\MSO \MCJ_F\ell_F\left(\GGW\right)\right)_\psi \right]'_\MSO & = & \prod_{\chi\in\MCX_\psi} \left[ \left(\MSO \ell_F\left(\GGW\right)\right)_\chi : \left(\MSO\MCJ_F\ell_F\left(\GGW\right)\right)_\chi \right]'_\MSO \nonumber \\
 & = & \prod_{\chi\in\MCX_\psi} \left[ \left(\MSO w_\infty L_{\MFf_\chi}\left(0,\bar{\chi}_{pr}\right)\right)_\chi : \left(\MSO\MCJ_F w_\infty L_{\MFf_\chi}\left(0,\bar{\chi}_{pr}\right)\right)_\chi \right]'_\MSO \nonumber \\
 & = & \prod_{\chi\in\MCX_\psi} \left[ \MSO [\mathrm{G}]_\chi : \left(\MSO\MCJ_F\right)_\chi \right]'_\MSO.
\label{equadiffmaisnon}
\end{eqnarray}
By Proposition \ref{prdirectsums} and Proposition \ref{propextscal pour modinc}, we have
\begin{equation}
\prod_{\chi\in\MCX_\psi} \left[ \MSO [\mathrm{G}]_\chi : \left(\MSO\MCJ_F\right)_\chi \right]'_\MSO = \left[\MSO[\mathrm{G}]_\psi : (\MSO\MCJ_F)_\psi\right]'_\MSO = \MSO\left[\MBZ_{\la g\ra}[\mathrm{G}]_\psi : \left(\MBZ_{\la g\ra}\MCJ_F\right)_\psi\right]'_{\MBZ_{\la g\ra}}
\label{aquaqua}
\end{equation}
By Corollary \ref{indmodinc}, we have
\begin{equation}
\left[\MBZ_{\la g\ra}[\mathrm{G}]_\psi : \left(\MBZ_{\la g\ra}\MCJ_F\right)_\psi\right]'_{\MBZ_{\la g\ra}} = \MBZ_{\la g\ra}\left[\MBZ_{\la g\ra}[\mathrm{G}]_\psi : \left(\MBZ_{\la g\ra}\MCJ_F\right)_\psi\right].
\label{equaqua}
\end{equation}
From (\ref{equadiffmaisnon}), (\ref{aquaqua}), (\ref{equaqua}) and the definition of $\MCJ_F$, we obtain
\[\left[ \left(\MSO \ell_F\left(\GGW\right)\right)_\psi : \left(\MSO \MCJ_F\ell_F\left(\GGW\right)\right)_\psi \right]'_\MSO = \MSO\left[\MBZ_{\la g\ra}[\mathrm{G}]_\psi : \left(\MBZ_{\la g\ra}\MCJ_F\right)_\psi\right] = \MSO\#\left(\MBZ_{\la g\ra}\otimes_\MBZ\mu(F)\right)_\psi.\]
Whence the corollary.
\hfill $\square$

To go further we need some preliminary remarks.

\begin{rmq}\label{rmqPsiHsubgrupG}
Let $H$ be a sub-group of $\mathrm{G}$.
Let $M$ and $N$ be two $\mathrm{G}$-modules, and $\psi:M\rightarrow N$ be a $\mathrm{G}$-equivariant map.
If $\mathrm{Coker}(\Psi):=N/\mathrm{Im}(\Psi)$ is annihilated by $\#(H)$ then we derive from $\Psi$ a surjective map
\[\Psi_\MSO:\MSO\otimes_\MBZ M\twoheadrightarrow\MSO\otimes_\MBZ N.\]

Let us assume, in addition, that $\mathrm{Ker}(\Psi)$ is annihilated by $\Sigma\GGs$, $\GGs\in H$.
Then, for every $\chi\in\widehat{\mathrm{G}}$ trivial on $H$, the restriction of $\Psi_\MSO$ gives an isomorphism
\[\left(\MSO\otimes_\MBZ M\right)_\chi \simeq \left(\MSO\otimes_\MBZ N\right)_\chi.\] 
\end{rmq}

As a particular case, for any subextension $K/k$ of $F/k$ and $H=\mathrm{Gal}(F/K)$, we shall consider the norm maps $Cl(\MCO_F)\rightarrow Cl(\MCO_K)$, $\mu(F)\rightarrow\mu(K)$, and also the map
\[\MBZ[\mathrm{G}]_0/\ell_F\left(\MCO_F^\times\right)\rightarrow\MBZ\left[\mathrm{Gal}(K/k)\right]_0/\ell_K\left(\MCO_K^\times\right),\]
deduced from the natural map $\MBZ[\mathrm{G}]\rightarrow\MBZ\left[\mathrm{Gal}(K/k)\right]$, where $\MBZ[\mathrm{G}]_0$ (resp. $\MBZ\left[\mathrm{Gal}(K/k)\right]_0$) is the augmentation ideal of $\MBZ[\mathrm{G}]$ (resp. $\MBZ\left[\mathrm{Gal}(K/k)\right]$).

\begin{rmq}\label{fittcardMSO}
For any commutative rings $A\subseteq B$ and any finitely generated $A$-module $M$, $Fitt_B(B\otimes_AM)=B.Fitt_A(M)$.
In particular, let $M$ be a finite $\MBZ[\mathrm{G}]$-module, let $\psi$ be an irreducible rational character of $\mathrm{G}$.
Then 
\[\MSO\#\left(\left(\MBZ_{\la g\ra}\otimes_{\MBZ}M\right)_\psi\right) = \prod_{\chi\in\MCX_\psi}  Fitt_\MSO\left(\left(\MSO\otimes_\MBZ M\right)_\chi\right).\]
\end{rmq}

\begin{pr}\label{indiceUOmega}
Let $\psi\neq1$ be an irreducible rational character of $\mathrm{G}$.
Then
\begin{equation}
\left[\left(\MSO\ell_F(\MCO_F^\times)\right)_\psi : \left(\MSO\ell_F(\GGW)\right)_\psi\right]'_{\MSO} = \MSO \frac {w_\infty^{dim(\psi)}\#\left(\MBZ_{\la g\ra}\otimes_\MBZ Cl(\MCO_F)\right)_\psi} {\#\left(\MBZ_{\la g\ra}\otimes_\MBZ\mu(F)\right)_\psi},
\label{equUOmega}
\end{equation}
where $\left(\MSO\ell_F(\MCO_F^\times)\right)_\psi$ and $\left(\MSO\ell_F(\GGW)\right)_\psi$ are viewed as $\MSO$-lattices of the $\MBC$-vector space $\MBC[\mathrm{G}]$.
\end{pr}

\noindent\textsl{Proof.}
Obviously, $F_\chi$ does not depend on the choice of $\chi\in\MCX_\psi$.
We set $F_\psi := F_\chi$, for $\chi\in\MCX_\psi$.
Let $\Xi_\psi$ be the set of $\chi\in\widehat{\mathrm{G}}$ such that $\mathrm{Ker}\,\chi$ strictly contains $\mathrm{Gal}\left(F/F_\psi\right)$.
For any $I\subseteq\Xi_\psi$, we define
\[F_I:=\lA\begin{array}{ll}
F_\psi & \text{if $\quad I=\vid$.} \\
\bigcap_{\chi\in I}F_\chi & \text{if $\quad I\neq\vid$.}
\end{array}\right.\]
Using Proposition \ref{prdirectsums} and the formula (\ref{equaMSOellGGWechi}), we have
\[\left[\MSO[\mathrm{G}]_\psi:\MSO\ell_F(\Omega)_\psi\right]'_\MSO = \prod_{\chi\in\MCX_\psi} \left[\MSO[\mathrm{G}]_\chi:\left(\MSO\ell_F(\GGW)\right)_\chi\right]'_\MSO = \prod_{\chi\in\MCX_\psi} \left[\MSO[\mathrm{G}]_\chi:\left(\MSO w_\infty L_{\MFf_\chi}(0,\bar{\chi}_{pr})\right)_\chi\right]'_\MSO.\]
By Proposition \ref{modincpourreseauxisomorphes}, and since $\#\left(\MCX_\psi\right)=dim(\psi)$,
\begin{equation}
\left[\MSO[\mathrm{G}]_\psi:\MSO\ell_F(\Omega)_\psi\right]'_\MSO = \MSO w_\infty^{dim(\psi)}\prod_{\chi\in\MCX_\psi}L_{\MFf_\chi}(0,\bar{\chi}_{pr}).
\label{equaprodLchu}
\end{equation}
For $I\subseteq\Xi_\psi$ fixed we denote by $\Xi_I$ the set of $\chi\in\widehat{\mathrm{G}}$ such that $F_\chi\subseteq F_I$.
The inclusion-exclusion principle and the analytic class number formula give
\begin{equation}
\prod_{\chi\in\MCX_\psi} L_{\MFf_\chi}(0,\bar{\chi}_{pr}) = \prod_{I\subseteq\Xi_\psi} \left(res(\zeta_k,0) \prod_{\substack{\chi\in\Xi_I \\ \chi\neq1}} L_{\MFf_\chi}(0,\bar{\chi}_{pr})\right)^{(-1)^{\#(I)}} = \prod_{I\subseteq\Xi_\psi} \left(res(\zeta_{F_I},0)\right)^{(-1)^{\#(I)}},
\label{equaprodresidu}
\end{equation}
where $\zeta_{F_I}$ (resp. $\zeta_k$) is the Dedekind zeta function of $F_I$ (resp. $k$) and $res(\zeta_{F_I},0)$ (resp. $res(\zeta_k,0)$) is the residue of $\zeta_{F_I}$ (resp. $\zeta_k$) at $0$.
Let $\MBZ\left[\mathrm{Gal}(F_I/k)\right]_0$ be the augmentation ideal of $\MBZ\left[\mathrm{Gal}(F_I/k)\right]$.
Then, it can be shown that
\begin{equation}
res(\zeta_{F_I},0) = - \frac{h\left(\MCO_{F_I}\right) \MCR\left(\MCO_{F_I}\right)}{w_{F_I}ln\left(q^d\right)},
\label{equaresidu}
\end{equation}
where $\MCR\left(\MCO_{F_I}\right):=\left[\MBZ\left[\mathrm{Gal}(F_I/k)\right]_0:\ell_{F_I}\left(\MCO_{F_I}^\times\right)\right]$. 
Since $\MBZ h\left(\MCO_{F_I}\right) = Fitt_\MBZ\left(Cl\left(\MCO_{F_I}\right)\right)$ we have
\begin{eqnarray*}
\MSO h\left(\MCO_{F_I}\right) = \MSO Fitt_\MBZ\left(Cl\left(\MCO_{F_I}\right)\right) & = &  Fitt_\MSO\left(\MSO\otimes_\MBZ Cl\left(\MCO_{F_I}\right)\right) \\
 & = & \prod_{\chi\in\Xi_I} Fitt_\MSO\left(\left(\MSO\otimes_\MBZ Cl\left(\MCO_{F_I}\right)\right)_\chi\right) \\ 
 & = & \prod_{\chi\in\Xi_I} Fitt_\MSO\left(\left(\MSO\otimes_\MBZ Cl\left(\MCO_F\right)\right)_\chi\right).
\end{eqnarray*}
The last equality is an application of Remark \ref{rmqPsiHsubgrupG}.
In the same manner we have 
\[\MSO w_{F_I} = \prod_{\chi\in\Xi_I} Fitt_\MSO\left(\left(\MSO\otimes_\MBZ \mu(F)\right)_\chi\right)\]
\[\text{and}\quad \MSO\MCR\left(\MCO_{F_I}\right) = \prod_{\chi\in\Xi_I} Fitt_\MSO\left(\left(\MSO\otimes_\MBZ \left( \MBZ[\mathrm{G}]_0 / \ell_F\left(\MCO_F^\times\right) \right)\right)_\chi\right).\]
For any $\chi\in\widehat{\mathrm{G}}$, let us set $h_\chi:=Fitt_\MSO\left(\left(\MSO\otimes_\MBZ Cl\left(\MCO_F\right)\right)_\chi\right)$, $w_\chi := Fitt_\MSO\left(\left(\MSO\otimes_\MBZ \mu(F)\right)_\chi\right)$, and $\MCR_\chi := Fitt_\MSO\left(\left(\MSO\otimes_\MBZ \left( \MBZ[\mathrm{G}]_0 / \ell_F\left(\MCO_F^\times\right) \right)\right)_\chi\right).$
Combining (\ref{equaprodresidu}) and (\ref{equaresidu}), and applying the inclusion-exclusion principle a second time we obtain
\begin{eqnarray}
\MSO\prod_{\chi\in\MCX_\psi} L_{\MFf_\chi}(0,\bar{\chi}_{pr}) = \MSO \prod_{I\subseteq\Xi_\psi} \left(\frac{h\left(\MCO_{F_I}\right)\MCR\left(\MCO_{F_I}\right)}{-w_{F_I}ln\left(q^d\right)}\right)^{(-1)^{\#(I)}} & = & \prod_{I\subseteq\Xi_\psi} \left(\prod_{\chi\in\Xi_I}h_\chi\MCR_\chi w_\chi^{-1}\right)^{(-1)^{\#(I)}} \nonumber \\
  & = & \prod_{\chi\in\MCX_\psi}h_\chi\MCR_\chi w_\chi^{-1}.
\label{equahchiRchiwchi}
\end{eqnarray}
By (\ref{equaprodLchu}), (\ref{equahchiRchiwchi}), Remark \ref{fittcardMSO}, Corollary \ref{indmodinc}, and Proposition \ref{propextscal pour modinc}, we have
\[\left[\MSO[\mathrm{G}]_\psi:\left(\MSO\ell_F(\Omega)\right)_\psi\right]'_\MSO = \frac{w_\infty^{dim(\psi)}\# \left( \MBZ_{\la g\ra} \otimes_\MBZ Cl\left(\MCO_F\right) \right)_\psi } {\# \left( \MBZ_{\la g\ra} \otimes_\MBZ \mu\left(F\right) \right)_\psi }\left[\MSO[\mathrm{G}]_\psi : \left(\MSO\ell\left(\MCO_F^\times\right)\right)_\psi\right]'_\MSO.\]
Multiplying by $\left[\left(\MSO\ell\left(\MCO_F^\times\right)\right)_\psi : \MSO[\mathrm{G}]_\psi\right]'_\MSO$, and applying Proposition \ref{propprodmodinc}, we have
\[\left[\left(\MSO\ell_F(\MCO_F^\times)\right)_\psi : \left(\MSO\ell_F(\GGW)\right)_\psi\right]'_{\MSO} = \MSO \frac {w_\infty^{dim(\psi)} \#\left(\MBZ_{\la g\ra}\otimes_\MBZ Cl(\MCO_F)\right)_\psi} {\#\left(\MBZ_{\la g\ra}\otimes_\MBZ\mu(F)\right)_\psi}.\]
\hfill $\square$

\begin{cor}\label{corOwmoinsunellGGW}
Let $\psi\neq1$ be an irreducible rational character of $\mathrm{G}$.
Then
\[\left[\left(\MSO\ell_F(\MCO_F^\times)\right)_\psi : \left(\MSO w_\infty^{-1}\ell_F(\GGW)\right)_\psi\right]'_{\MSO} = \MSO \frac {\#\left(\MBZ_{\la g\ra}\otimes_\MBZ Cl(\MCO_F)\right)_\psi} {\#\left(\MBZ_{\la g\ra}\otimes_\MBZ\mu(F)\right)_\psi}.\]
\end{cor}

\noindent\textsl{Proof.}
We have the equality
\[\left[ \left(\MSO \ell_F\left(\MCO_F^\times\right)\right)_\psi : \left(\MSO w_\infty^{-1}\ell_F\left(\GGW\right)\right)_\psi \right]'_\MSO = \]
\[\left[ \left(\MSO \ell_F\left(\MCO_F^\times\right)\right)_\psi : \left(\MSO \ell_F\left(\GGW\right)\right)_\psi \right]'_\MSO \left[\left(\MSO \ell_F\left(\GGW\right)\right)_\psi : \left(\MSO w_\infty^{-1}\ell_F\left(\GGW\right)\right)_\psi \right]'_\MSO,\]
thanks to Proposition \ref{propprodmodinc}.
By Proposition \ref{modincpourreseauxisomorphes}, (iii), we have
\[\left[\left(\MSO \ell_F\left(\GGW\right)\right)_\psi : \left(\MSO w_\infty^{-1}\ell_F\left(\GGW\right)\right)_\psi \right]'_\MSO = w_\infty^{-dim(\psi)}.\]
We conclude by using Proposition \ref{propprodmodinc} and the above Proposition \ref{indiceUOmega}.
\hfill $\square$

\begin{theo}\label{theocalcind}
Let $\psi\neq1$ be an irreducible rational character of $\mathrm{G}$.
We have
\[\left[\left(\MBZ_{\la g\ra}\otimes_\MBZ\MCO_F^\times\right)_\psi : \left(\MBZ_{\la g\ra}\otimes_\MBZ\MCE_F\right)_\psi\right] = \#\left(\MBZ_{\la g\ra}\otimes_\MBZ Cl\left(\MCO_F\right)\right)_\psi.\]
\end{theo}

\noindent\textsl{Proof.}
We deduce from the identity $\MCO_F^\times\cap \mathrm{Ker}(\ell_F)=\mu(F)$ that the factor $\mathrm{G}$-modules $\MCO_F^\times/\MCE_F$ and $\ell_F\left(\MCO_F^\times\right)/\ell_F\left(\MCE_F\right)$ are isomorphic.
In particular,
\begin{equation}
\left[\left(\MBZ_{\la g\ra}\otimes_\MBZ\MCO_F^\times\right)_\psi : \left(\MBZ_{\la g\ra}\otimes_\MBZ\MCE_F\right)_\psi\right] = \left[\left(\MBZ_{\la g\ra} \ell_F\left(\MCO_F^\times\right)\right)_\psi : \left(\MBZ_{\la g\ra} \ell_F\left(\MCE_F\right)\right)_\psi\right].
\label{equasansellell}
\end{equation}
Let us also remark that, by the property (ii) of Stark units, $\MCP_F^{\GGs-1}\subseteq\MCE_F$ for all $\GGs\in \mathrm{G}$.
Thus, since $\psi\neq1$ and $g\in\MSO^\times$ we obtain
$\left(\MSO \ell_F\left(\MCE_F\right)\right)_\psi = \left(\MSO \ell_F\left(\MCP_F\right)\right)_\psi$.
But $\ell_F(\MCP_F)$ and $\ell_F(\GGW)$ are related by the equality
$w_\infty\ell_F(\MCP_F) = \MCJ_F\ell_F(\GGW)$,
so that
\begin{equation}
\left(\MSO \ell_F\left(\MCE_F\right)\right)_\psi = \left(\MSO w_\infty^{-1}\MCJ_F\ell_F\left(\GGW\right)\right)_\psi.
\label{equaMCEGGW}
\end{equation}
Therefore, Corollary \ref{indmodinc}, Proposition \ref{propextscal pour modinc} and the formula (\ref{equaMCEGGW}) give
\[\begin{array}{lll}
\MSO\left[\left(\MBZ_{\la g\ra}\ell_F\left(\MCO_F^\times\right)\right)_\psi : \left(\MBZ_{\la g\ra}\ell_F\left(\MCE_F\right)\right)_\psi\right] & = & \MSO\left[\left(\MBZ_{\la g\ra} \ell_F\left(\MCO_F^\times\right)\right)_\psi : \left(\MBZ_{\la g\ra} \ell_F\left(\MCE_F\right)\right)_\psi\right]'_{\MBZ_{\la g\ra}}\\
 & = & \left[\left(\MSO\ell_F\left(\MCO_F^\times\right)\right)_\psi : \left(\MSO w_\infty^{-1}\MCJ_F\ell_F(\GGW)\right)_\psi\right]'_\MSO.
\end{array}\]
By Proposition \ref{propprodmodinc}, this last $\MSO$-index-module is the product of the two index-modules already computed in Corollary \ref{corOMegaJFOmega} (see also Remark \ref{trucidiot}) and Corollary \ref{corOwmoinsunellGGW}.
Thus we obtain 
\[\MSO\left[\left(\MBZ_{\la g\ra}\otimes_\MBZ\MCO_F^\times\right)_\psi : \left(\MBZ_{\la g\ra}\otimes_\MBZ\MCE_F\right)_\psi\right] = \MSO\#\left(\left(\MBZ_{\la g\ra}\otimes_\MBZ Cl\left(\MCO_F\right)\right)_\psi\right).\]
Since the integers we are comparing are prime to $g$, the theorem follows.
\hfill $\square$
\bigskip

\noindent\textbf{Acknowledgement.}
I wish to express all my gratitude to Hassan Oukhaba, who introduced me to this topic.

\bibliographystyle{amsplain}

%% Fin %%%%%%%%%%%%%%%%%%%%%%%%%%%%%%%%%%%%%%%%%%%%%%%%%%%%%
\end{document}